\documentclass[12pt,reqno]{article}

\usepackage[usenames]{color}
\usepackage{amssymb}
\usepackage{amsmath}
\usepackage{amsthm}
\usepackage{enumitem}
\usepackage{amsfonts}
\usepackage{amscd}
\usepackage{graphicx}
\usepackage{float}

\usepackage{url}
\usepackage[T1]{fontenc}
\usepackage[utf8]{inputenc}

\usepackage[colorlinks=true,
linkcolor=webgreen,
filecolor=webbrown,
citecolor=webgreen]{hyperref}

\definecolor{webgreen}{rgb}{0,.5,0}
\definecolor{webbrown}{rgb}{.6,0,0}

\usepackage{color}
\usepackage{fullpage}
\usepackage{float}

\usepackage{graphics}
\usepackage{latexsym}
\usepackage{epsf}
\usepackage{breakurl}

\newcommand{\seqnum}[1]{\href{https://oeis.org/#1}{\rm \underline{#1}}}

\begin{document}

\theoremstyle{plain}
\newtheorem{theorem}{Theorem}
\newtheorem{corollary}[theorem]{Corollary}
\newtheorem{lemma}[theorem]{Lemma}
\newtheorem{proposition}[theorem]{Proposition}

\theoremstyle{definition}
\newtheorem{definition}[theorem]{Definition}
\newtheorem{example}[theorem]{Example}
\newtheorem{conjecture}[theorem]{Conjecture}

\theoremstyle{remark}
\newtheorem{remark}[theorem]{Remark}

\begin{center}
\vskip 1cm{\Large\bf 
Gaps Between Consecutive Primes and the\\
\vskip .1in
Exponential Distribution
}

\vskip 1cm
Joel E. Cohen\\
The Rockefeller University\footnote{Author's other affiliations:  Departments of Statistics, Columbia University \& University of Chicago.}\\
1230 York Avenue, Box 20\\
New York, NY 10065\\
USA\\
\href{mailto: cohen@rockefeller.edu }{\tt cohen@rockefeller.edu }\\
ORCID 0000-0002-9746-6725
\ \\
\end{center}

\vskip .2 in

\def\modd#1 #2{#1\ \mbox{\rm (mod}\ #2\mbox{\rm )}}

\vskip .2in
\begin{abstract}
Based on the primes less than $4\times 10^{18}$,
Oliveira e Silva {et al. (\emph{Math. Comp.},
83(288):2033–2060, 2014)} conjectured
an asymptotic formula for the sum of the $k$th power 
of the gaps between consecutive primes less than a large number $x$.
We show that the conjecture of Oliveira e Silva holds if and only if
the $k$th moment of the first $n$ gaps 
is asymptotic to the $k$th moment
of an exponential distribution with mean $\log n$,
though the distribution of gaps is not exponential.
Asymptotically exponential moments imply that the gaps asymptotically 
obey Taylor's law of fluctuation scaling: 
variance of the first $n$ gaps $\sim$ (mean of the first $n$ gaps)$^2$.
If the distribution of the first $n$ gaps is asymptotically exponential with mean $\log n$,
then the expectation of the largest of the first $n$ gaps is asymptotic to $(\log n)^2$.
The largest of the first $n$ gaps is asymptotic to $(\log n)^2$ 
if and only if the Cram\'er-Shanks conjecture holds.
Numerical counts of gaps and the maximal gap $G_n$ among 
the first $n$ gaps test these results.
While most values of $G_n$ are better approximated by 
$(\log n)^2$ than by other models, 
seven values of $n$ with $G_{n} >2e^{-\gamma}(\log n)^2$
suggest that
{$\limsup_{n \to\infty} G_n/[2e^{-\gamma}(\log n)^2]$ may exceed 1.}
\end{abstract}


\noindent \emph{Disclosure}: The author reports there are no competing interests to declare.
\section{Introduction}
The gaps between consecutive prime numbers have been studied mathematically, numerically, and statistically
by many people for centuries
\cite{dickson1919, HL23, cramer1936, shanks1964, gallagher1976, MAIER1981257,
 heath-brown_1982, Maier1985221, guy1994, RieselHans1994Pnac,yamasaki1994,granville1995,wolf1998,cutter2001,CP2005,wolf2014, oliveira2014,cohen2016,wolf2017, enwiki:1167645599, kourbatov2019upper}. 
This list of references is far from complete.
Many questions remain open.
Throughout, we use \emph{gap} to mean a {difference, denoted $d$,}
between \emph{consecutive} primes, including $d_4 := p_5 - p_4 = 11-7=4$
and excluding $p_4 - p_2 = 7-3=4$,
{where $p_n$ is the $n$th prime, $p_1 = 2$.}

We do not attempt 
to summarize here everything that is currently known about the distribution of prime gaps. 
{Some selected recent results and conjectures are pertinent.}
Heath-Brown \cite[p.\ 87]{heath-brown_1982} conjectured a simple asymptotic expression for the sum
of the squares of the gaps of the primes not exceeding a large real $x$ (OEIS \seqnum{A074741}).
Wolf \cite[Conjecture 5]{wolf1998} conjectured an alternative asymptotic form and compared the two alternatives numerically.
Based on a detailed numerical study of the primes less than $4\times 10^{18}$,
Oliveira e Silva \cite[p.\ 2056]{oliveira2014} conjectured a simple asymptotic formula for
the sum of the $k$th power
of the gaps of primes less than $x$ as $x\to\infty$,
for every $k=1, 2,\ldots$.
Wolf \cite{wolf2017} refined and extended the conjectured asymptotic expressions of Oliveira e Silva
and compared the alternatives analytically and numerically.
{Firoozbakht conjectured that the sequence 
$(p^{1/n}_n)_{n\in \mathbb{N}}$ decreases as $n$ increases \cite{enwiki:1167645599}.
The conjecture holds for all primes less than $2^{64}\approx 1.84 \times 10^{19}$ as of April 2024 \cite{enwiki:1167645599}. 
Kourbatov \cite{kourbatov2019upper} showed that Firoozbakht's conjecture 
implies an upper bound on all gaps.
Maier \cite{MAIER1981257} showed that chains of consecutive large gaps exist, 
and that average gaps of primes in short intervals do not converge
(contrary to previous conjecture) 
but have a limsup strictly greater than the liminf \cite{Maier1985221}.
}

What is new here? We consider the first $n=1, 2,\ldots$ gaps 
instead of the gaps of primes less than large $x$ (section \ref{sec:exponential}).
This change of scale enables us to show (Theorem \ref{th:moments}) 
that the conjecture of Oliveira e Silva holds if and only if,
for every $k=1, 2,\ldots$,
the $k$th moment of the gaps 
(that is, the sum of the $k$th power of the first $n$ gaps, divided by $n$)
is asymptotically (as $n\to\infty$) 
the $k$th moment 
$k!(\log n)^k$
of an exponential distribution with mean $\log n$.

Obviously the distribution of gaps is not exponential because
an exponential random variable has a continuous density function on $[0,\infty)$,
whereas gaps are positive integers only 
{and all gaps except $d_1 = p_2 - p_1 = 3 - 2 = 1$ are even integers}.
Gaps are also not geometrically distributed (the integer-valued analog of the exponential)
because a geometric distribution would give positive probability to all odd 
{positive integers.}
Despite these deviations from the exponential distribution at a microscopic scale,
the numerically supported conjecture that the $k$th moment of the first $n$ gaps 
is asymptotically the $k$th moment $k!(\log n)^k$
of an exponential distribution with mean $\log n$ has interesting consequences, 
as we now summarize and then show in detail.

From Theorem \ref{th:moments}, it follows (Corollary \ref{th:TL}) that, 
as $n\to\infty$, the ratio of the variance of the first $n$ gaps 
to the square of the mean of the first $n$ gaps converges to 1 (section \ref{sec:Taylor}).
Thus the gaps have a power-law asymptotic variance function
\cite{tweedie1946,tweedie1947,tweedie1984,barlevenis1986, barlevstramer1987, davidiancarroll1987} 
or, equivalently, {asymptotically} obey Taylor's law of fluctuation scaling:
{variance of the first $n$ gaps $\sim$ 
square of the mean of the first $n$ gaps}
\cite{taylor2019taylor,eisler2008fluctuation,cohen2016,cohen2023,cohenmichael2017,demers2018}.
{Moreover (section \ref{sec:maxgap}),
Taylor's law (with a different coefficient and a different exponent)
also describes the variance function of the maximal gap among the first $n$ gaps
if the maximal gap may be modeled
as the largest order statistic of $n$ observations of an exponential distribution
with mean $\log n$: variance of the maximal order statistic $\sim (\pi^2/6) \times$mean of the maximal order statistic.
These (conjectured) scaling laws do} not appear to have been noticed previously.

The Cram\'er-Shanks conjecture 
{(Cram\'er \cite[pp.\ 24, 27]{cramer1936}, Shanks \cite[p.\ 648, his Eq. (5)]{shanks1964}, 
Granville \cite[p.\ 21, his Eq. (14)]{granville1995})}
proposes that the maximal gap of
primes less than $x$ is asymptotic to $(\log x)^2$ as $x\to\infty$ (section \ref{sec:maxgap}). 
Using the exponential distribution as a heuristic model of the gaps,
we show (Theorem \ref{th:expasymptotic}) that the Cram\'er-Shanks conjecture holds if and only if
the largest of the first $n$ gaps is asymptotic (as $n\to\infty$) to
the expectation $(\log n)^2$
of the largest order statistic of a sample of $n$ independent observations 
from an exponential distribution with mean $\log n$.
{However,} the $(1-1/n)$-quantile of the largest order statistic of a sample of $n$ independent observations 
from an exponential distribution with mean $\log n$
yields $2(\log n)^2$ as
an asymptotic estimate of the largest of the first $n$ gaps.
{This estimate exceeds the Cram\'er-Shanks conjecture and Kourbatov's
bound \cite{kourbatov2019upper} on gaps derived from Firoozbakht's conjecture.}

In the concluding section \ref{sec:numerical},
we compare asymptotic expressions 
(derived from the exponential model of prime gaps and from other conjectures)
for moments and maximal gaps 
with exact numerical results derived from public sources.
Numerical evidence supports the moments of an exponential distribution with mean $\log n$
as models of the moments of the first $n$ gaps
and questions $(\log n)^2$ as an asymptotic description of
the largest of the first $n$ gaps.

\section{Definitions and background}

For any natural number $n\in\mathbb{N}:=\{1, 2, \ldots \}$,
let $p_n$ be the $n$th prime starting from $p_1=2,\ p_2=3$ (OEIS \seqnum{A000040}). 
Let $\mathbb{P}:=\{p_1, p_2,  \ldots \}$ be the set of primes.
Define the $n$th gap $d_n$ between consecutive primes (OEIS \seqnum{A001223}) 
and the 
sum $D_k(x)$ of the $k$th power, $k\in \{0\}\cup \mathbb{N}$, of the gaps 
between consecutive primes that do not exceed real $x$ to be
\begin{align}
d_n&:=p_{n+1}-p_n,\ n\in\mathbb{N},\label{eq:defgap}\\
D_k(x)&:=\sum_{p_{n+1\leq x}}d_n^k.\label{eq:Dkx}
\end{align}
{We shall say that $d_n$ is the gap of $p_n$ if and only if \eqref{eq:defgap}.
For example, the gap of 17 is 2.

Let $\pi(x):= \#\{p\in\mathbb{P} \mid p \leq x \}$ be the number of primes that do not exceed real $x\geq 2$ 
(OEIS \seqnum{A000720} {gives $\pi(x)$ evaluated at positive} integral values of $x$).
Obviously $\pi(p_n)=n$ and $p_{\pi(x)}= \max\{p_{n}\in\mathbb{P} \mid p_{n}\leq x\}$. 
For $2<a<b$, the number of primes in the interval $(a,b]$
equals $\pi(b) - \pi(a)$.
Those primes are $p_{\pi(a)+1}, \ldots,p_{\pi(b)}$
and the gaps of those primes are $d_{\pi(a)+1}, \ldots,d_{\pi(b)}$.
Then 
\begin{align}\label{eq:sumgaps}
d_{\pi(a)+1} + \cdots + d_{\pi(b)-1} := L_1 < L:=b - a < L_2 :=d_{\pi(a)+1} + \cdots + d_{\pi(b)}.
\end{align}
Some authors define the number of the gaps of the primes in $(a,b]$ to be
$N_1 := \pi(b)-\pi(a) - 1 $ 
(excluding gap $d_{\pi(b)}$ because $p_{\pi(b)+1}>b$)
and some authors define it to be $N_2 := \pi(b)-\pi(a)$
(including gap $d_{\pi(b)})$.
Asymptotically, the difference in these definitions is immaterial.
}

If $f(x)$ and $g(x)$ are real-valued functions of real $x$ and $g(x)>0$ for all 
$x$ sufficiently large, define $f(x)\sim g(x)$ to mean that 
$f(x)/g(x) \to 1$ as $x\to\infty$.
In Riesel's \cite[p.\ 61]{RieselHans1994Pnac} cautious notation, define $f(x) \sim_c  g(x)$ if $f(x)\sim g(x)$ is conjectured but not proved.
Likewise, $=_c$ will denote a conjectured but unproved equality.

Oliveira e Silva \cite[p.\ 2056]{oliveira2014} conjectured that
\begin{align}
D_k(x)\sim_c k!\ x (\log x)^{k-1},\qquad k \geq 1. \label{eq:OeSconjecture}
\end{align}

The prime number theorem states:
\begin{align}
\pi(x)\sim \frac{x}{\log x}\text{ as }\ x\to\infty \text{ or equivalently }p_n \sim n \log n\text{ as }n\to\infty.\label{eq:PNT}
\end{align}
Because $\pi(x)\sim {x}/{\log x}$
if and only if $\log x \sim x/\pi(x)$, and because $x/\pi(x)$ is asymptotic to the 
average gap of the primes $\leq x$, 
the prime number theorem \eqref{eq:PNT} implies that the average gap of the primes $\leq x$
is asymptotic to $\log x$.
{
(Why? Because $x$ is asymptotic to the sum of gaps of the primes $\leq x$.
Also $\pi(x)$ is asymptotic 
to the number of gaps of the primes $\leq x$,
by either definition of that number $N_1$ or $N_2$ above.
By definition, the average of the gaps of the primes $\leq x$ is 
the sum of gaps of the primes $\leq x$ divided by 
the number of gaps of the primes $\leq x$.
That ratio is asymptotic to $ x/\pi(x) = \log (x)$.)
}
As $n\to\infty$, \eqref{eq:PNT} implies that $p_n/n\to\infty$.
{Taking logs of both sides gives} $\log p_n - \log n \to\infty$ but 
$\log p_n / \log n \ {\sim (\log n + \log \log n)/\log n} \to 1$.

Wolf \cite{wolf2017} remarked that, when $k=0$, the {right side of} conjecture \eqref{eq:OeSconjecture}  
{becomes $x/\log x$ as in}
the prime number theorem in \eqref{eq:PNT}.
When $k=1$, \eqref{eq:Dkx} simplifies to 
$D_1(x) = p_{\pi(x)}-2$ and
\eqref{eq:OeSconjecture} asserts that 
$D_1(x) \sim x$,
which is well known and, in any case,  
follows from Lemma \ref{le:primeratio} below.
For $k>1$, \eqref{eq:OeSconjecture} apparently remains unproved, though well supported numerically
\cite{oliveira2014, wolf2017}.

We shall use some elementary consequences of the prime number theorem \eqref{eq:PNT}.

\begin{lemma} \label{le:primeratio}
The prime number theorem \eqref{eq:PNT} implies
\begin{align}
\lim_{n\to\infty} \frac{p_{n+1}}{p_n}&=1,\label{eq:ratsucprimes}\\
\lim_{x\to\infty}\frac{x}{p_{\pi(x)} }&=1, \label{eq:ratxtop}\\
\lim_{x\to\infty}\frac{\log \pi(x)}{\log x}&=1. \label{eq:logxsimlogn}
\end{align}
\end{lemma}
\begin{proof}
By \eqref{eq:PNT},
\begin{align}
\frac{p_{n+1}}{p_n}&\sim \frac{(n+1)\log(n+1)}{n\log n}\sim 1\cdot 1=1.
\end{align}
This proves \eqref{eq:ratsucprimes}.
To prove \eqref{eq:ratxtop}, observe that 
$p_{\pi(x)}\leq x < p_{\pi(x)+1}$, 
hence $1\leq x/p_{\pi(x)} < p_{\pi(x)+1}/p_{\pi(x)}$.
As $x\to\infty$, also $\pi(x)\to\infty$, so \eqref{eq:ratsucprimes} implies \eqref{eq:ratxtop}.
To prove \eqref{eq:logxsimlogn}, take logarithms of \eqref{eq:PNT}:
$\log \pi(x) \sim \log x - \log\log x \sim \log x$ as $x\to\infty$.
\end{proof}

Now we define an exponentially distributed random variable 
and state, {mostly} without proof, several of its well known properties.
Let $X$ be a nonnegative real-valued random variable with cumulative distribution function {(cdf)}
$F(x):=\Pr\{X \leq x\},\ x \geq 0$.
If there exists some $\lambda>0$ such that $F(x)=1-\exp(-\lambda x)$, 
then $X$ is \emph{exponentially distributed} with scale parameter $\lambda$ 
and we write $X\stackrel{d}{=}Exp(\lambda)$.

\begin{lemma} \label{le:exp}
Assume $\lambda>0$. {For $n>1$, let $X$ and $X_1,\ldots,X_{n}$ be independently and identically distributed (iid)
$Exp(\lambda)$ random variables. Then:}
\begin{enumerate}
\item \label{le:exp1} For any real $c>0$, $cX\stackrel{d}{=}Exp(\lambda/c)$.
\item \label{le:moments} For real $r>-1$,
the $r$th moment of $X$ is $E(X^r)=\Gamma(r+1)/\lambda^r$ \cite[p.\ 28, Eq. (5)]{marshallolkin2007}.
When $r\in\mathbb{N}$, then $\Gamma(r+1)=r!$.
{Hence the expectation or mean of $X$ is $1/\lambda$ 
and the variance of $X$ is $1/\lambda^2$.}
\item \label{le:exp3} Let the order statistics of 
$X_1,\ldots,X_{n}$ be $X_{(1)}\leq \cdots \leq X_{(n)}$.
Then \cite[p.\ 343]{sarhangreenberg1962} for $i=1,\ldots,n$,
\begin{align}
E(X_{(i)})&=\frac{1}{\lambda}\sum_{j=1}^i \frac{1}{n-j+1}=\frac{1}{\lambda}\sum_{j=n-i+1}^n \frac{1}{j},\label{eq:mn}\\
Var(X_{(i)})&=\frac{1}{\lambda^2}\sum_{j=1}^i \frac{1}{(n-j+1)^2}=\frac{1}{\lambda^2}\sum_{j=n-i+1}^n \frac{1}{j^2}.\label{eq:vr}
\end{align}

\item \label{le:expleast} The cumulative distribution function (cdf) of the 
{smallest order statistic $X_{(1)}$ is
\begin{align}
\Pr\{X_{(1)}\leq x\}=1-\exp(-n\lambda x),\ x\geq 0,
\end{align}
i.e., $X_{(1)}\stackrel{d}{=}Exp(\lambda n)$.
For $0<q<1$, let $y_1$ be the $(1-q)$-quantile of $X_{(1)}$. 
By definition, $y_1$ satisfies 
$q = \Pr\{X_{(1)}>y_1\}=1-\Pr\{X_{(1)}\leq y_1\}=\exp(-n\lambda y_1) $.
Then
\begin{align}\label{eq:qqy}
y_1 = - (1/[n\lambda])\log q < + (1/[n\lambda])(1/q-1).
\end{align}
In particular, if $q = 1/n$, then $y_1 = + (1/[n\lambda])\log n$.

\item \label{le:exp4}
The cdf of the} largest order statistic $X_{(n)}$ is
\begin{align}
\Pr\{X_{(n)}\leq x\}=(1-\exp(-\lambda x))^n,\ x\geq 0.
\end{align}
{For $0<q<1$, let $y_n$ be the $(1-q)$-quantile of $X_{(n)}$. 
By definition, $y_n$ satisfies $\Pr\{X_{(n)}>y_n\}=1-\Pr\{X_{(n)}\leq y_n\}=q$.
Then
\begin{align}\label{eq:qqy}
y_n = - (1/\lambda)\log(1 - (1-q)^{1/n}) < +(1/\lambda)( \log n - \log q).
\end{align}
\begin{proof}
$q = 1 - (1-\exp(-\lambda y_n))^n$, hence
$1-\exp(-\lambda y_n) = (1-q)^{1/n}$ and then
$-\lambda y_n = \log(1 - (1-q)^{1/n})$.
The inequality follows from $(1-q)^{1/n}<1 - q/n$.
\end{proof}
}
\item \label{le:exp5} 
\cite[p.\ 324, Ex. 9]{grimmettstirzaker2001}
{By Lemma \ref{le:exp}.1 above, $\lambda X_n \stackrel{d}{=}Exp(1)$. Then}
\begin{align}
\Pr\left\{\limsup_{n\to\infty} \frac{\lambda X_n}{\log n}=1\right\}=1.
\end{align}

\item \label{le:exp2} 
Let $S_n:=X_1+\cdots+X_{n}$. 
Let random variables $U_1,\ldots,U_{n-1}$ be iid uniformly on $[0,1]$ and 
let $U_0:=0,\ U_n:=1$ with probability 1. Then 
$(X_1/S_n,X_2/S_n,\ldots,X_{n}/S_n)\stackrel{d}{=}
(g_1, g_2, \ldots, g_{n}):=(U_{(1)}-U_{(0)}, U_{(2)}-U_{(1)},\ldots,U_{(n)}-U_{(n-1)})$ 
(\cite[p.\ 75, III.3(e)]{feller1971} and \cite[p.\ 302, Ex. 42]{grimmettstirzaker2001}).

{
\item \label{le:largedev} \cite{MITlargedeviations}
Large deviations: When $a>1/\lambda$, then
\begin{align}
\Pr\left\{\frac{X_1+\cdots+X_{n}}{n} > a\right\}\approx \exp(-(a\lambda - 1 - \log\lambda-\log a)n).
\end{align}
}
\end{enumerate}
\end{lemma}

\section{Asymptotic moments of first $n$ gaps are exponential}\label{sec:exponential}
Define the $k$th moment ($k>-1$) of the first $n$ gaps to be
\begin{align}
\mu_{k,n}':=\frac{1}{n}\sum_{j=1}^n d_j^k.\label{eq:kthmomngaps}
\end{align}
Thus if $\pi(x)=n$, then \eqref{eq:Dkx} and \eqref{eq:kthmomngaps} give $\mu_{k,n}'=D_k(x)/n$.
\begin{theorem}\label{th:moments}
The conjecture \eqref{eq:OeSconjecture} of Oliveira e Silva \cite[p.\ 2056]{oliveira2014}
holds for each $k\in \{0\}\cup \mathbb{N}$ if and only if
the $k$th moment of the first $n$ gaps is asymptotic to
the $k$th moment of the exponential distribution $Exp(1/\log n)$ with mean $\log n$:
\begin{align}
\mu_{k,n}' \sim_c k! (\log n)^k,\qquad n\to\infty.                   \label{eq:kthmomngapsExp}
\end{align}
\end{theorem}

\begin{proof}
The cases $k=0,\ k=1$ are known to be true. 
Assume \eqref{eq:OeSconjecture} for $k>1$. 
Using the definition \eqref{eq:Dkx} and Lemma \ref{le:primeratio},
replace $x$ by $p_{n+1}$ in \eqref{eq:OeSconjecture} to get
\begin{align}\label{eq:keyproof}
D_k(p_{n+1})&\sim_c k!\ n \log n (\log(n \log n))^{k-1}=k!\ n \log n (\log n + \log\log n)^{k-1}\sim k!\ n (\log n)^k.
\end{align}
Dividing both sides by $n$ gives \eqref{eq:kthmomngapsExp}.

Conversely, suppose the $k$th moment of the first $n$ gaps is asymptotic to
the $k$th moment of the exponential distribution $Exp(1/\log n)$.
Then, working backward through the calculation in \eqref{eq:keyproof},
\begin{align}
n\mu_{k,n}'&:=\sum_{j=1}^n d_j^k \sim k!n(\log n)^k
\sim k!\ n \log n (\log n + \log\log n)^{k-1} \nonumber\\
&=k!\ n \log n (\log(n \log n))^{k-1}
\sim_c D_k(p_{n+1})
.\label{eq:converse1}
\end{align}
Again using Lemma \ref{le:primeratio} to replace $p_{n+1}$ by $x$ and $n\log n \sim p_n \sim p_{n+1}$ by $x$ gives
\begin{align}
k!\ x (\log x)^{k-1}
\sim_c D_k(x)
.  \label{eq:converse2}
\end{align}
\end{proof}

\section{Asymptotic variance function obeys Taylor's power law}\label{sec:Taylor}
For $n\in\mathbb{N},\ n\geq 2$, define the mean and variance of the first $n$ gaps as
\begin{align}
m_n&:=\frac{1}{n}\sum_{j=1}^n d_j=\mu_{1,n}', \label{eq:mean} \\ 
v_n&:=\frac{1}{n-1}\sum_{j=1}^n (d_j-m_n)^2
=\frac{n}{n-1}\left(\frac{1}{n}\sum_{j=1}^n d_j^2-m_n^2\right)=\frac{n}{n-1}(\mu_{2,n}'-(\mu_{1,n}')^2).\label{eq:varexpanded}
\end{align}
{The central equality in \eqref{eq:varexpanded} holds because
$\sum_{j=1}^n (d_j-m_n)^2 
= \sum_{j=1}^n (d_j^2-2d_jm_n + m_n^2)
= \sum_{j=1}^n d_j^2 - 2 \sum_{j=1}^nd_jm_n+nm_n^2
= \sum_{j=1}^n d_j^2 - 2 (nm_n)m_n+nm_n^2
=  \sum_{j=1}^n d_j^2 -nm_n^2.$
}

The variance function of gaps is the function $f$ that satisfies $v_n=f(m_n)$ 
\cite{tweedie1946,tweedie1947,tweedie1984,barlevenis1986, barlevstramer1987, davidiancarroll1987}.
The asymptotic variance function of gaps is a (not unique) function $f$ that satisfies $v_n\sim f(m_n)$.

\begin{corollary}\label{th:TL}
Conditional on the conjectures of Heath-Brown \cite{heath-brown_1982}, Wolf \cite{wolf1998},
and \eqref{eq:OeSconjecture} of Oliveira e Silva \cite[p.\ 2056]{oliveira2014},
\begin{align}
\lim_{n\to\infty}\frac{v_n}{m_n^2}=_c1.\label{eq:TLgaps}
\end{align}
\end{corollary}
\begin{proof}
From \eqref{eq:kthmomngapsExp},
\begin{align}
\frac{\mu_{2,n}'}{(\mu_{1,n}')^2}\sim_c 2.
\end{align}
Hence
\begin{align}\label{eq:asymptTL}
v_n\sim (\mu_{1,n}')^2\left(\frac{\mu_{2,n}'}{(\mu_{1,n}')^2}-1\right)
\sim_c (\mu_{1,n}')^2(2-1)=m_n^2.
\end{align}
\end{proof}

This power-law asymptotic variance function is the special case $c=1,\ b=2$ 
of Taylor's law of fluctuation scaling $v_n\sim cm_n^b,\ c>0$
\cite{taylor2019taylor,eisler2008fluctuation,cohen2016,cohen2023,cohenmichael2017,demers2018, cohen2020brokenstickTL}.
These parameter values $c=1,\ b=2$ are consistent with 
the exponential distribution, 
in which the variance equals the square of the mean.
However,
a variance asymptotic to the squared mean does not imply an exponential distribution 
(e.g., a normal distribution with mean $\mu$ and variance $\mu^2$ has a variance equal to the squared mean).

\section{An asymptotically exponential distribution of gaps is plausible}
In a heuristic and numerical exploration, Yamasaki and Yamasaki \cite{yamasaki1994}
``assume that the exponential distribution
can be applied to the gaps of prime numbers. But the gaps [except for $d_1 = 1$] are always even
integers, so do not distribute continuously on $[0, \infty)$.''
They proposed adjustments to the exponential distribution.

Conditional on a uniform version of the unproved prime $r$-tuple conjecture of 
Hardy and Littlewood \cite{HL23},
Gallagher \cite{gallagher1976} showed that,
for $\lambda>0,\ h\sim \lambda \log N$, and $n\leq N$, 
as $N\to\infty$,
the distribution of the values of 
$\pi(n+h)-\pi(n)$ converges to the Poisson distribution with parameter $\lambda$.
Thus $\lambda$ is the expectation and the variance of the number of primes over the interval 
{$(n, n + \lambda\log N]$}
of length $\lambda \log N$, 
so the average gap is asymptotically $\log N$.
If the Poisson-distributed number of primes arrived in the interval $(n, n + \lambda\log N]$
according to a Poisson process,
then the inter-arrival interval between {successive} primes (that is, the gap) would be exponentially distributed
\cite[pp.\ 188, 378]{feller1971}, $Exp(1/\log N)$.

Wolf \cite{wolf2014} gave strong numerical evidence \cite[p.\ 3, his Figure 1]{wolf2014} and
a heuristic argument \cite[p.\ 2, his Eq. (7)]{wolf2014}
that an exponential distribution with mean $\log n$ asymptotically approximates the distribution of the first $n$ gaps.
To circumvent the difficulty that the exponential distribution lacks the discreteness of the natural numbers,
Wolf \cite[his Eq. (19), his Fig. 5]{wolf2014} proposed a rescaling that turns the discrete $d_n$ into a continuous variable.

Our finding that,
conditional on the conjecture \eqref{eq:OeSconjecture} of Oliveira e Silva \cite[p.\ 2056]{oliveira2014},
the moments of the first $n$ gaps are asymptotic to the moments of the exponential distribution,
even though the gaps are not continuously distributed,
is consistent with T. Stieltjes's discovery in 1894 that the moments of the lognormal distribution
do not uniquely specify the continuous probability distribution that produced the moments
\cite[pp.\ 17--18, section 2.1]{schmudgen2020}.

We suggest it is plausible that the moments \eqref{eq:kthmomngaps}
of the first $n$ gaps are asymptotic to the moments of $Exp(1/\log n)$.
First, the primes that do not exceed $x$ are asymptotically 
uniformly distributed on $[0,x]$ as $x\to\infty$ \cite{cohen2016}.
To see why, choose any $0<r<1$. 
Then, for large $x$, the number of primes not exceeding $rx$ divided by the number of primes not exceeding $x$
is asymptotically, by the prime number theorem \eqref{eq:PNT}, $\pi(rx)/\pi(x)\sim\lim_{x\to\infty} (rx/\log(rx))/(x/\log x) = r$.
(From a more general perspective, the asymptotic counting function $\pi(x)\sim x/\log x$ 
of primes is regularly varying \cite{bingham1987} with exponent 1.)

Assume $U_1,\ldots,U_{n}$ are iid uniform on $[0,x]$.
The finding that the primes are asymptotically uniformly distributed on $[0,x]$ \cite{cohen2016}
suggests the order statistics 
$U_{(1)}\leq \cdots \leq U_{(n)}$
of $U_1,\ldots,U_{n}$
as a stochastic model of the asymptotic distribution of the primes on $[0,x]$
and suggests the $n$ ``spacings'' \cite{Holst1980,cohen2020brokenstickTL}
$g_i:=U_{(i)}-U_{(i-1)},\ i=1,\ldots,n$, with $U_{(0)}:=0$ almost surely,
as a stochastic model of the first $n$ gaps.
The spacings are exponentially distributed 
(Lemma \ref{le:exp}.\ref{le:exp2})
and their average size is asymptotically $\log n$.
So the gaps between consecutive primes $\leq x$ are plausibly distributed as the 
spacings between consecutive order statistics of $\pi(x)$ iid uniform random variables on $[0,x]$, 
and these spacings are distributed as $Exp(1/\log n)$.

This stochastic model of the primes and gaps intentionally omits the discreteness of integers
and the lack of independence of primes.

\section{The largest of the first $n$ gaps}\label{sec:maxgap}

What would the exponential model of gaps predict about the asymptotic behavior of the maximal gap,
$G_n:=\max\{d_1,\ldots,d_n\}$, among the first $n$ gaps? 
Here we consider two possible answers.

\begin{theorem}\label{th:expasymptotic}
{For fixed $n>1$,}
let $X_i\stackrel{d}{=}Exp(1/\log n),\ i=1,\ldots,n$, be iid.
Let $\gamma\approx 0.5772$ be the Euler-Mascheroni constant,
$\zeta(\cdot)$ the zeta function, and $\pi^2/6\approx 1.6449$.
Then the mean and variance of the largest order statistic $X_{(n)}$ are
\begin{align}
E(X_{(n)})&=(\log n)\sum_{j=1}^n \frac{1}{n-j+1}\sim (\log n)(\gamma+\log n)\sim (\log n)^2,\label{eq:largestmean}\\
Var(X_{(n)})&=(\log n)^2\sum_{j=1}^n \frac{1}{(n-j+1)^2}\sim (\log n)^2\zeta(2)= (\log n)^2\frac{\pi^2}{6}.\label{eq:largestvar}
\end{align}
Thus $Var(X_{(n)})/E(X_{(n)})\to \pi^2/6>1$ and $Var(X_{(n)})\sim (\pi^2/6)E(X_{(n)})$.
The maximal gap illustrates Taylor's law with $c=\pi^2/6,\ b=1$.
\end{theorem}
\begin{proof}
From Lemma \ref{le:exp}.\ref{le:exp3}, 
\begin{align}
E(X_{(n)})&=(\log n)\left(1+\frac{1}{2} +\frac{1}{3} +\cdots+\frac{1}{n} \right)
\sim (\log n)(\gamma+\log n)\sim (\log n)^2,\label{eq:largestmean2}\\
Var(X_{(n)})&=(\log n)^2\left(1+\frac{1}{2^2} +\frac{1}{3^2} +\cdots+\frac{1}{n^2} \right)
\sim (\log n)^2\zeta(2)= (\log n)^2\frac{\pi^2}{6}.\label{eq:largestvar2}
\end{align}
\end{proof}

The Cram\'er-Shanks conjecture may be stated as
\begin{align}
\max_{p_j \leq x} (p_{j+1}-p_j)\sim_c (\log x)^2
\sim (\log p_{\pi(x)})^2.\label{eq:cramershanksconjecture}
\end{align}
See Cram\'er \cite[pp.\ 24, 27]{cramer1936}, Shanks \cite[p.\ 648, his Eq. (5)]{shanks1964}, 
and the very helpful Granville \cite[p.\ 21, his Eq. (14)]{granville1995}.

\begin{theorem}\label{th:Cramerorderstats}
The Cram\'er-Shanks conjecture \eqref{eq:cramershanksconjecture} holds if and only if
the maximal gap $G_n:=\max\{d_1,\ldots,d_n\}$ satisfies 
\begin{align}
G_n\sim_c (\log n)^2.\label{eq:lognsquared}
\end{align}
\end{theorem}
\begin{proof}
If there are $n$ primes less than or equal to $x$, 
then the largest of these is $p_n\sim n\log n \sim x$, using Lemma \ref{le:primeratio}.
Thus, translated from the scale of primes not exceeding $x$
to the first $n$ gaps,
the Cram\'er-Shanks conjecture \eqref{eq:cramershanksconjecture} becomes \eqref{eq:lognsquared} because
\begin{align}
\max_{1\leq j \leq n} d_j\sim_c (\log (n \log n))^2=(\log n + \log\log n)^2\sim(\log n)^2.
\end{align}

Conversely, if $\max_{1\leq j \leq n} d_j \sim_c (\log n)^2$, then 
$\max_{1\leq j \leq n} d_j=\max_{p_j \leq p_n} d_j\sim\max_{p_j \leq x} d_j$
and $\log x \sim \log (n\log n) \sim \log n + \log\log n \sim \log n$,
so $(\log x)^2 \sim(\log n)^2$ and \eqref{eq:cramershanksconjecture} holds.
\end{proof}

Theorems \ref{th:expasymptotic}
and \ref{th:Cramerorderstats}
suggest that $(\log x)^2$ in the Cram\'er-Shanks conjecture \eqref{eq:cramershanksconjecture}
is too small asymptotically.
In the model of Theorem \ref{th:expasymptotic},  
the largest gap is likely to be larger than the expectation of the largest gap 
because the variance of the largest gap is of the same order of magnitude as the expectation of the largest gap.
To the extent that the above conjectures are correct that gaps have exponential moments asymptotically,
it would be surprising if the largest of the first $n$ gaps did not exceed $(\log x)^2$.

{These suggestions indirectly challenge Firoozbakht's conjecture \cite{enwiki:1167645599}.
Kourbatov \cite{kourbatov2019upper} showed that Firoozbakht's conjecture implies that 
$d_n:=p_{n+1} - p_n <_c (\log p_n)^2 - \log p_n - 1$ for all $n > 9$
and is implied by
$d_n <_c (\log p_n)^2 - \log p_n - 1.17$ for all $n > 9$.
These upper bounds are stronger than the Cram\'er-Shanks conjecture
\eqref{eq:cramershanksconjecture}, and lower than the dominant term
$(\log p_n)^2\sim (\log n)^2$ (see Theorem \ref{th:Cramerorderstats})
by amounts that increase without limit as $n$ increases.
To the extent that our Theorems \ref{th:expasymptotic}
and \ref{th:Cramerorderstats} 
suggest that $(\log x)^2$ in the Cram\'er-Shanks conjecture \eqref{eq:cramershanksconjecture}
is too small asymptotically,
our results suggest even more that Kourbatov's upper bounds on gaps are likely to be too small.
Kourbatov \cite{kourbatov2019upper} discusses many related inequalities and conjectures concerning prime gaps.
}

Granville \cite[p.\ 24]{granville1995} reviewed results that suggest,
contrary to the Cram\'er-Shanks conjecture \eqref{eq:cramershanksconjecture}, that 
\begin{align}\label{eq:granvilleconjecture}
\max_{p_j \leq x} d_j\stackrel{>}{\sim}_c 2e^{-\gamma}(\log x)^2
\sim 2e^{-\gamma}(\log p_{\pi(x)})^2, 
\end{align}
where $2e^{-\gamma}\approx 1.1229$.
Translating from the scale of $x$ to the scale of the number of gaps $n=\pi(x)$ by $x\sim p_n \sim n\log n$ so that
$(\log x)^2\sim (\log(n \log n))^2\sim (\log n)^2$
gives the proposal that
\begin{align}
G_n\stackrel{>}{\sim}_c 2e^{-\gamma}(\log n)^2.\label{eq:Granvillen}
\end{align}

These observations motivate another estimate of the largest of the first $n$ gaps,
namely, the $(1-1/n)$-quantile of the largest order statistic $X_{(n)}$.
We now determine the value of $y$ such that,
in a sample of size $n$ from $Exp(1/\log n)$,
as in Theorem \ref{th:expasymptotic},  
the probability that the largest order statistic $X_{(n)}$ exceeds $y$ is $1/n$.

\begin{theorem}\label{th:uptail2}
{In a sample of size $k > 1$ from $Exp(1/\log n)$,
let $y_k$ be the $(1-1/k)$-quantile of the maximal order statistic $X_{(k)}$.
Then $y_k \sim 2(\log n)(\log k)$ as $k\to\infty$. 
If $k \sim tn$ for some $t>0$, then 
$y_k\sim2(\log n)^2$ as $n\to\infty$.
}
\end{theorem}

\begin{proof}
Using Lemma \ref{le:exp}.\ref{le:exp4}, \eqref{eq:qqy}, with $\lambda=1/\log n$ gives 
\begin{align}
{y_k = -(\log n)\log(1-(1-1/k)^{1/k})}.
\end{align}
{Repeated applications of l'Hopital's rule as $k\to\infty$ shows that
$-\log(1-(1-1/k)^{1/k}) \sim (2 \log k)$
and hence that
\begin{align}
y_k \sim +2(\log n)(\log k). \label{eq:x}
\end{align}
}
\end{proof}

Wolf \cite[his Eq. (15)]{wolf2014} argued heuristically for the conjecture that, as $x\to\infty$,
\begin{align}
G_W(x):=\max_{p_n<x}(p_n-p_{n-1})\sim_c \frac{x}{\pi(x)}(2\log(\pi(x))-\log(x)+c),\label{eq:WolfGx}
\end{align}
where $c \approx 0.2778769$.
{Here $c := \log(C_2)$ and $C_2 := 2 \prod_{p>2}(1- (p-1)^{-2})$ 
is known sometimes as the twins constant 
and sometimes as twice the twins constant.
The product in $C_2$ is taken over all positive odd primes.}
Replacing $x$ by $p_n+\varepsilon$ and letting $\varepsilon \downarrow 0$ gives the equivalent conjecture, as $n\to\infty$,
\begin{align}
G_W(p_n)\sim_c \frac{p_n}{n}(2\log(n)-\log(n \log n)+c)\sim (\log n)^2.\label{eq:WolfGpn}
\end{align}

\section{Numerical illustrations}\label{sec:numerical}
For moments and maximal gaps, we compare exact numerical results 
derived from public sources with asymptotic expressions 
derived from the exponential model of prime gaps and from other conjectures.

\subsection{Moments}

To illustrate Theorem \ref{th:moments},
we compare exact computations of the moments of the first $n$ gaps with the 
moments of an exponential distribution with mean $\log n$ as given by \eqref{eq:kthmomngapsExp}.
Wolf \cite{wolf2014} generously made public extensive numerical results
at \url{http://pracownicy.uksw.edu.pl/mwolf/gapstau.zip}.
Wolf \cite[p.\ 2, his Eq. (2)]{wolf2014} defined, for large real $x$ and \emph{even} $d\in\mathbb{N}$,
\begin{align}
\tau_d(x):=\#\{p_{n+1} \mid p_{n+1}<x\quad \mathrm{ and }\quad p_{n+1}-p_n=d \}
\end{align}
and, for every \emph{odd} $d\in\mathbb{N}$, $\tau_d (x):=0$.
In words, when $d$ is even, $\tau_d(x)$ is the number of pairs of consecutive primes 
such that the larger prime is less than $x$
and such that their difference equals $d$;
and for odd $d$, $\tau_d (x):=0$.
This definition and my subsequent analyses exclude 
$\tau_1 (x)=1,\ x\geq 3$, representing the single gap $p_2-p_1=3-2$ of size $d=1$.
Wolf collected the non-zero values of $\tau_d(x)$ in 34 text files,
one file for each of the 34 values of $x=2^t,\ t=15(1)48$.
(Approximately, $2^{48}\approx 2.8147 \times 10^{14}$.)
Each file has two columns, the first listing values of $d$ such that $\tau_d (x)>0$,
and the second listing the corresponding non-zero values of $\tau_d (x)$.

First, using MATLAB Version 9.13.0.2049777 (R2022b), 
I calculated $\tau_d(x)$ directly for all primes less than $2^{20}$. 
My results agreed exactly with those in Wolf's file \texttt{tau20s.dat}. 
For example, Wolf and I independently found  $\tau_2(2^{20})=8535$ gaps of size $d=2$;
the maximal gap size, $d=114$, occurred once: $\tau_{114}(2^{20})=1$.

\begin{figure}[H]
\centering
\includegraphics[width=\textwidth]{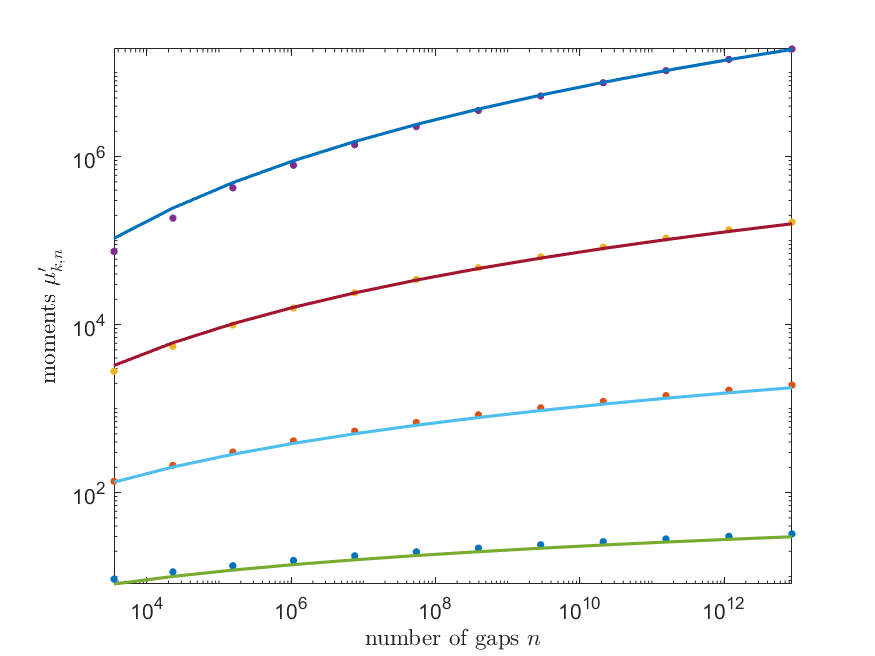}
\caption{Moments $\mu_{k,n}'\ (k=1, 2, 3, 4)$ of the first $n$ gaps between consecutive primes (solid dots)
from $k=1$ (bottom row) to $k=4$ (top row);
and corresponding $k$th moments $(k=1, 2, 3, 4)$ of exponential distributions $Exp(1/\log n)$ (solid lines).
The dots and the lines are calculated independently with no 
adjustment of parameters.
}
\label{fig:1}
\end{figure}

Having partially verified Wolf's results, I analyzed 
Wolf's 12 files with values of $x=2^t,\ t=15, 18, 21, 24, 27,\ldots,48$.
For each such file separately, 
Table \ref{tab:1} gives $t$; the number of gaps
$n=\sum_{d\in\mathbb{N}}\tau_d(x)$ where $x=2^t$ (not counting the first gap $d=1$);
the first four integer moments $\mu_{k,n}',\ k=1, 2, 3, 4$ from \eqref{eq:kthmomngaps};
and the maximal gap $G_n$.
Figure \ref{fig:1} compares 
$\mu_{k,n}'$ with $k! (\log n)^k$,
which is asymptotic to 
the corresponding moments of samples of size $n$ from an exponential distribution $Exp(1/\log n)$,
as conjectured in \eqref{eq:kthmomngapsExp}.
Visually, the agreement in Figure \ref{fig:1} between the exact moments and the asymptotic moments of 
$Exp(1/\log n)$ is remarkable.

\subsection{Maximal gaps}
We compare the first 80 known maximal gaps (now \emph{including} the initial gap of size 1)
\cite{enwiki:primegap}
with four asymptotic models: 
\eqref{eq:lognsquared}, \eqref{eq:cramershanksconjecture},
\eqref{eq:Granvillen}, and \eqref{eq:granvilleconjecture}.
These results on maximal gaps are more extensive and more detailed
than those in Table \ref{tab:1}.
The tabulation lists three sequences:
\begin{enumerate}
\item the indices of primes after which a maximal gap occurs \\
$n=1,2,4,9,24,30,99,154,189,217,1183,\ldots$ (OEIS \seqnum{A005669}),
\item the corresponding maximal gaps \\
$G_n:=p_{n+1}-p_n=1,2,4,6,8,14,18,20,22,34,36,\ldots$ (OEIS \seqnum{A005250}), and
\item the primes after which a maximal gap occurs \\
$p_n=2, 3, 7, 23, 89, 113, 523, 887, 1129, 1327, 9551,\ldots$ (OEIS \seqnum{A002386}).
\end{enumerate}

I partially verified the Wikipedia tabulation in two ways. First,
Marek Wolf (personal communication, 2022-11-10) generously sent me his tabulation
of the first 74 maximal gaps $G_n$ (not including $G_1$), their indices $n$, and their beginning primes $p_n$.
Second, using MATLAB, I calculated these three sequences directly for all primes less than $10^6$. 
My results and Wolf's agreed exactly, as far as they went, with the Wikipedia tabulation \cite{enwiki:primegap}. 

\begin{figure}[H]
\centering
\includegraphics[width=\textwidth]{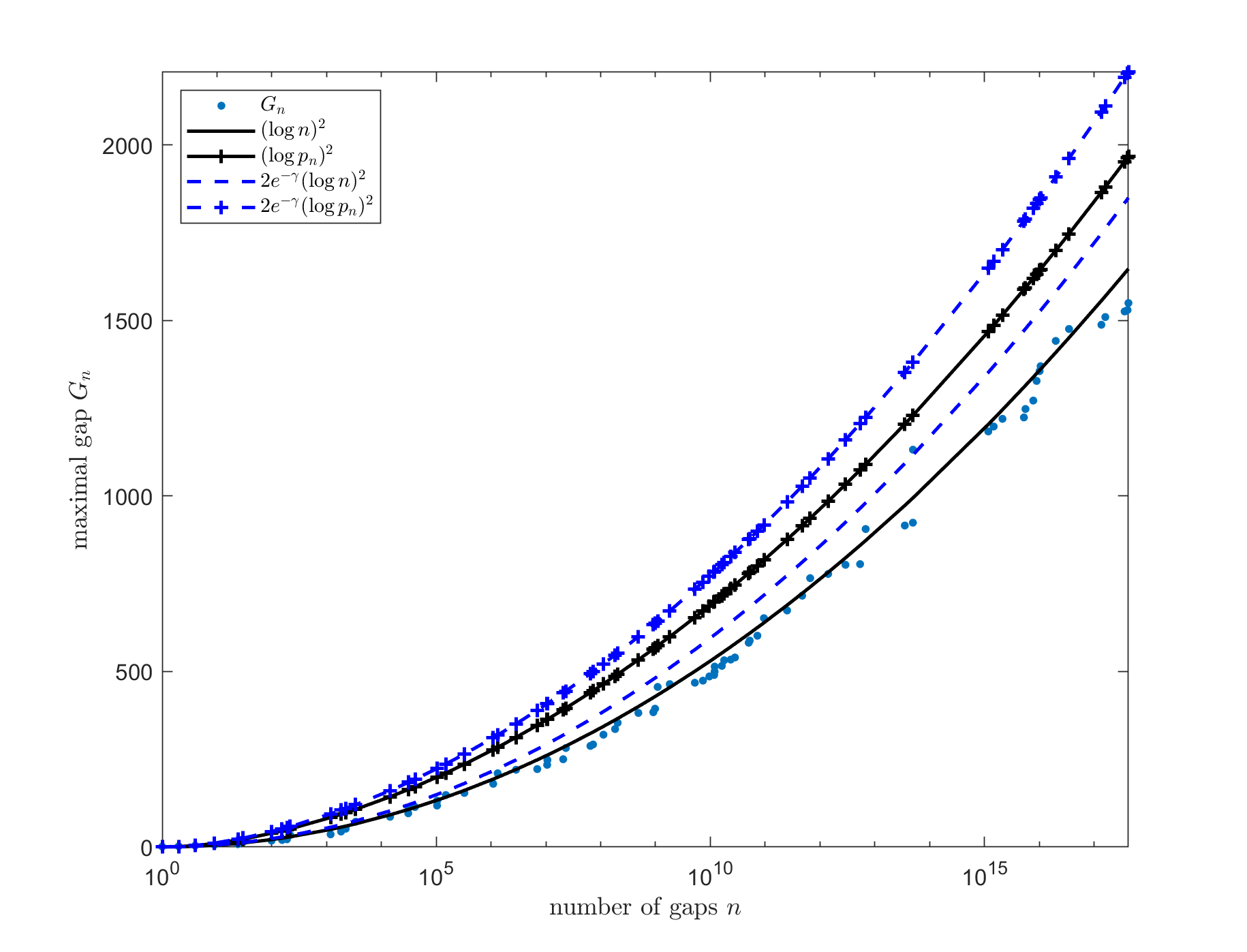}
\caption{Maximal gap $G_n$ (solid blue dots) 
among the first $n$ gaps {and four conjectured models:} 
$(\log n)^2$ (solid black line);
$(\log p_n)^2$ (solid black line with $+$ markers);
$2e^{-\gamma}(\log n)^2$ (dashed purple line), 
where $2e^{-\gamma}\approx 1.1229$;
and $2e^{-\gamma}(\log p_n)^2$ (dashed purple line with $+$ markers).
The dots and lines are calculated independently with no 
adjustment of parameters.
{While} the great majority of the values of $G_n$ are better approximated by 
$(\log n)^2$ (solid black line) {than by the other models}, 
{the seven values of $n$ with $G_{n} >2e^{-\gamma}(\log n)^2$
in Table \ref{tab:2} suggest that}
{$\limsup_{n \to\infty} G_n/[2e^{-\gamma}(\log n)^2]$ may exceed 1.}
}
\label{fig:2}
\end{figure}

Figure \ref{fig:2} plots, as a function of $n$,
the maximal gap size $G_n$ among the first $n$ gaps; 
and the values of the four conjectured asymptotic expressions
$(\log n)^2$;
$(\log p_n)^2$;
$2e^{-\gamma}(\log n)^2$, where $2e^{-\gamma}\approx 1.1229$;
and $2e^{-\gamma}(\log p_n)^2$.
The majority of points $(n, G_n)$ (blue solid dots in Figure \ref{fig:2}) fall 
slightly below the conjectured asymptotic behavior
$(n, (\log n)^2)$ (solid black line in Figure \ref{fig:2}) 
derived here from the conjectured asymptotically exponential distribution of gaps.
The form of Wolf's conjecture in \eqref{eq:WolfGpn}
falls slightly below $(n, (\log n)^2)$ and often 
closer to the points $(n, G_n)$ than $(n, (\log n)^2)$.
However, 
{the seven values of $n$ with $G_{n} >2e^{-\gamma}(\log n)^2$
in Table \ref{tab:2} suggest that}
{$\limsup_{n \to\infty} G_n/[2e^{-\gamma}(\log n)^2]$ may exceed 1.}
{The available numerical results neither confirm nor reject}
the suggestion implied by Theorem \ref{th:uptail2} that  
{$\limsup_{n \to\infty} G_n$}
might be asymptotic to $2(\log n)^2$.
The asymptotic behavior of $G_n$ remains to be determined.

\begin{table}
{
\begin{tabular} {|r|r|r|r r r r|r|}\hline
{ } & {} & {}&\multicolumn{4}{|c|}{{moments}} &{} \\ 
{row} & {$t$} & {$n$} & $k=1$ &$k=2$ &$k=3$ &$k=4$ & {$G_n$} \\ \hline
{1} & {15} & {3510} & {9.3293} & {136.2017} & {2.7818e+03} & {7.4292e+04} & {72} \\ 
{2} & {18} & {22998} & {11.3982} & {210.7095} & {5.5060e+03} & {1.8546e+05} & {86} \\ 
{3} & {21} & {155609} & {13.4770} & {304.1124} & {9.8914e+03} & {4.2503e+05} & {148} \\ 
{4} & {24} & {1077869} & {15.5652} & {412.7866} & {1.5776e+04} & {7.8862e+05} & {154} \\ \hline
{5} & {27} & {7603551} & {17.6520} & {539.4491} & {2.3885e+04} & {1.3864e+06} & {222} \\ 
{6} & {30} & {54400026} & {19.7379} & {683.2373} & {3.4423e+04} & {2.2806e+06} & {282} \\ 
{7} & {33} & {393615804} & {21.8231} & {844.1273} & {4.7670e+04} & {3.5440e+06} & {354} \\ 
{8} & {36} & {2.8744e+09} & {23.9074} & {1.0222e+03} & {6.3972e+04} & {5.2773e+06} & {464} \\ \hline
{9} & {39} & {2.1152e+10} & {25.9908} & {1.2173e+03} & {8.3638e+04} & {7.5819e+06} & {532} \\ 
{10} & {42} & {1.5666e+11} & {28.0736} & {1.4296e+03} & {1.0699e+05} & {1.0574e+07} & {652} \\ 
{11} & {45} & {1.1667e+12} & {30.1560} & {1.6590e+03} & {1.3435e+05} & {1.4377e+07} & {766} \\ 
{12} & {48} & {8.7312e+12} & {32.2379} & {1.9056e+03} & {1.6603e+05} & {1.9127e+07} & {906} \\ 
\hline
\end{tabular}
}
\caption{For each upper limit $x=2^t$, this table shows
the exponent $t$ of 2, the number $n$ of gaps between consecutive
primes (not counting the odd first gap), the moments $\mu_{k,n}' (k=1, 2, 3, 4)$, and the maximal gap $G_n$. 
}
\label{tab:1}
\end{table}

\begin{table}
{
\begin{tabular} {|r|r|r|r|r|r|r|r|}\hline
{row} & {$n$} & {$G_n$} & {$p_n$} & {$(\log n)^2$} & {$(\log p_n)^2$} & {$2e^{-\gamma}(\log n)^2$} & {$2e^{-\gamma}(\log p_n)^2$}  \\ \hline
{1} & {1             } & {1   } & {2               } & {0} & {0.4805} & {0} & {0.5395} \\ 
{2} & {2             } & {2   } & {3               } & {0.4805} & {1.2069} & {0.5395} & {1.3553} \\
{3} & {4             } & {4   } & {7               } & {1.9218} & {3.7866} & {2.1580} & {4.2520} \\ \hline
{4} & {9             } & {6   } & {23              } & {4.8278} & {9.8313} & {5.4212} & {11.0398} \\ 
{5} & {30            } & {14  } & {113             } & {11.5681} & {22.3482} & {12.9901} & {25.0952} \\ 
{6} & {217           } & {34  } & {1327            } & {28.9433} & {51.7058} & {32.5010} & {58.0614} \\ 
{7} & {49749629143526} & {1132} & {1693182318746371} & {994.6470} & {1229.6} & {1116.9} & {1380.7} \\ 
\hline
\end{tabular}
}
\caption{Seven maximal gaps $G_n$ that exceed $(\log n)^2$ and $2e^{-\gamma}(\log n)^2$, 
the index $n$ of the prime $p_n$ that begins the maximal gap, 
the prime $p_n$ that begins the maximal gap, 
$(\log p_n)^2$, and $2e^{-\gamma}(\log p_n)^2$. 
Only for $n=1$ and $n=2$ is $G_n>2e^{-\gamma}(\log p_n)^2$.
}
\label{tab:2}
\end{table}

\section*{Acknowledgements}
J.E.C. thanks Andrew Granville, D. R. Heath-Brown, Tom\'as Oliveira e Silva, and Marek Wolf for very helpful guidance to relevant references; Marek Wolf for sharing numerical results analyzed here; and Roseanne Benjamin for assistance during this work.
{Two helpful reviewers provided constructive suggestions.}

\bibliographystyle{plain}
\bibliography{refs20240416}

\bigskip
\hrule
\bigskip

\noindent (Concerned with sequences
\seqnum{A000040}, 
\seqnum{A000720},
\seqnum{A001223},
\seqnum{A074741}.)

\bigskip
\hrule
\bigskip

\end{document}